\documentclass[letterpaper, 10 pt, conference]{ieeeconf}  

\IEEEoverridecommandlockouts                              

\usepackage{epsfig} 

\usepackage{amsmath} 
\usepackage{amssymb}  
\usepackage{tabularx}
\usepackage{amsmath}
\usepackage{balance}

\newcounter{subassumption}[asu]

\makeatletter
\renewcommand{\p@subassumption}{\theasu}
\makeatother

\usepackage{rotating,booktabs,multirow}
\makeatletter
\let\NAT@parse\undefined
\makeatother
\usepackage{hyperref}
\usepackage{graphicx}
\usepackage{subfigure}
\newtheorem{remark}{Remark}[section]

\title{\LARGE \bf
Model-free immune therapy: A control approach to acute inflammation$^{*}$
}

\author{Ouassim Bara$^{1}$,  Michel Fliess$^{2, 5}$, C\'{e}dric Join$^{3, 5, 6}$, Judy Day$^4$, Seddik M. Djouadi$^{1}$ 
\thanks{*Work partially supported by the NSF-DMS Award 1122462  and a Fullbright Scholarship.}
\thanks{$^{1}$Department of Electrical Engineering  and Computer Science,
        University of Tennessee, Knoxville, TN 37996, USA. \newline
        {\tt \small \{obara, mdjouadi\}@utk.edu}}
\thanks{$^{2}$LIX (CNRS, UMR 7161), \'Ecole polytechnique, 91128 Palaiseau, France. {\tt \small Michel.Fliess@polytechnique.edu } } 
\thanks{$^3$CRAN (CNRS, UMR 7039)), Universit\'{e} de Lorraine, BP 239, 54506 Vand{\oe}uvre-l\`{e}s-Nancy, France. \newline
{\tt\small cedric.join@univ-lorraine.fr}}
\thanks{$^{4}$Department of Mathematics, University of Tennessee, Knoxville, TN 37996, USA.
        {\tt \small judyday@utk.edu}} 
\thanks{$^{5}$AL.I.E.N. (ALg\`{e}bre pour Identification \& Estimation Num\'{e}riques), 24-30 rue Lionnois, BP 60120, 54003 Nancy, France. \newline
        {\tt \small \{michel.fliess, cedric.join\}@alien-sas.com}}  
\thanks{$^{6}$Projet NON-A, INRIA Lille -- Nord-Europe, France. }    
}

\begin{document}
\maketitle
\thispagestyle{empty}
\pagestyle{empty}

\begin{abstract}
Control of an inflammatory immune response is still an ongoing research. Here, a strategy consisting of manipulating a pro and anti-inflammatory mediator is considered. Already existing and promising model-based techniques  suffer unfortunately from a most difficult calibration. This is due to the different types of inflammations and to the strong parameter variation between patients. This communication explores another route via the new model-free control and its corresponding 
``intelligent'' controllers. A ``virtual'' patient, \textit{i.e.}, a mathematical model, is only employed for digital simulations. A most interesting feature of our control strategy is the fact that the two outputs which must be driven are sensorless. This difficulty is overcome by assigning suitable reference trajectories to two other outputs with sensors. Several most encouraging computer simulations, corresponding to different drug treatment strategies, are displayed and discussed.

\end{abstract}

\section{Introduction}
The importance, complexity and ubiquity of the notions of \emph{infection} and \emph{inflammation} are well explained by the following quotation \cite{nathan}: 
\textit{The `inflammatory process' includes a tissue-based startle reaction to trauma; go/no-go decisions based on integration of molecular clues for tissue penetration by microbes; the beckoning, instruction and dispatch of cells; the killing of microbes and host cells they infect; liquefaction of surrounding tissue to prevent microbial metastasis; and the healing of tissues damaged by trauma or by the host's response. If at any step an order to proceed is issued but progress to the next step is blocked, the inflammatory process may detour into a holding pattern, such as infiltration of a tissue with aggregates of lymphocytes and leukocytes (granulomas) that are sometimes embedded in proliferating synovial fibroblasts (pannus), or distortion of a tissue with collagen bundles (fibrosis). Persistent inflammation can oxidize DNA badly enough to promote neoplastic transformation. } 
According to \cite{cohen}, \textit{the overall mortality is approximately $30\%$, rising to $40\%$ in the elderly and is $50\%$ or greater in patients with the more severe syndrome}. 
The corresponding literature is of course huge. See, \textit{e.g}, 
\begin{itemize}
\item \cite{schwartz}  on the cause, 
\item \cite{virchow,david,israel,schulz,virchow-or,warburg} for the connections with cancer,
\item \cite{bricaire,deeks,hunt,nowak} for the interactions with the human immunodeficiency virus (HIV). 
\item \cite{dantzer,raison} for the possible relationship with depression.
\end{itemize}
Although applying automatic control to immune therapy has attracted some interest, as depicted in \cite{parker}, it is much less developed than in other domains, like, \textit{e.g.}, for insulin-dependent diabetes (see, \textit{e.g.}, \cite{beq,doyle}, and the references therein).  Let us nevertheless mention promising papers using respectively optimal control (\cite{bara2015,acc16,kirschner1997,stengel,stengel1,stengel2}) and predictive control (\cite{day2010using,hogg, zitelli}).
Those approaches are model-based. Among those papers, the most recent ones (\cite{bara2015,acc16,day2010using,zitelli}) use the same set of phenomenological ordinary differential equations from \cite{reynolds2006} (see 
also \cite{reynolds} and \cite{day}):
\begin{itemize}
\item The corresponding model is based on the  non-specific protective mechanism, namely, the \emph{innate immune response}, in contrast to the \emph {adaptive immune system}. The latter provides a more advanced and strategic response producing $B$ and $T$ cells together with specific antibodies.\footnote{See, \textit{e.g.}, the classic textbook \cite{murphy} for an explanation of the technical medical words here, and elsewhere in this communication.}
\item Anti-inflammatory mediators are included. They play an important r\^{o}le to mitigate a severe inflammation and, therefore, avoid tissue damage and high pathogen proliferation.
\item Its biological relevance has been confirmed via a good qualitative reproduction of severe systemic inflammation in a biological organism. 
\end{itemize}
Other mathematical modelings have been proposed (see, \textit{e.g.}, \cite{al1993disease,esaim,ho,kumar,perelson,rundell,song,y}).
In spite of interesting preliminary results in \cite{cdc13, baraparameter,zitelli}, state observation and parameter identification are not yet fully mastered. Its calibration, which depends heavily on the type of inflammatory response and on patient differences (genetics, age, gender, \dots), is therefore 
most intricate. 

This paper suggests another route, namely the recent model-free setting and the corresponding ``intelligent'' controllers \cite{csm}.\footnote{See also \cite{automatica,gao,larminat}.}  It is worthwhile to recall that model-free control has already been successfully 
applied in quite diverse case-studies (see, \textit{e.g.}, \cite{toulon,berlin,siam,med16} in the field of ``life engineering''). The modeling remains nevertheless irreplaceable 
at this stage for \textit{in silico} testing, \textit{i.e.}, for computer simulations. We will also be employing \cite{reynolds2006}. Let us emphasize the following key point: there is no need for the proposed control
technique to use any state observer and any parameter identification technique. 

From a purely control-theoretic standpoint, a major novelty of this study lies in the necessity to drive sensorless states. The poor knowledge of the system makes the derivation of an observer quite intractable. The solution lies in a ``good understanding'' of the system, \textit{i.e.}, in the design of an ``efficient'' reference trajectories tracking with respect to the states with sensors. Such a feedforward ``philosophy'' is of course inspired by flatness-based control (see \cite{flmr}, and 
\cite{murray,levine,sira-flat}).


Our paper is organized as follows. Sections \ref{model} and \ref{mfc} review respectively the mathematical modeling and model-free control. Several computer simulations are displayed and discussed in Section \ref{virtual}. Suggestions for future research may be found in Section \ref{conclusion}.

\section{A virtual patient}\label{model}
A mathematical model, \textit{i.e.}, a \emph{virtual patient}, via four ordinary differential equations, for an acute inflammatory response to pathogenic infection has been proposed \cite{reynolds2006}:
{\small{\begin{align}
\hspace{-3em}\frac{dP}{dt} &= k_{pg} P (1-\frac{P}{P_{\infty}})-\frac{k_{pm} s_{m}P}{\mu_{m}+k_{mp} P}-k_{pn}f(N)P\label{eq:mod1}\\  
\hspace{-3em}\frac{dN}{dt} &= \frac{s_{nr}R}{\mu_{nr}+R} - \mu_n N+ u_p(t) \label{eq:mod2}\\ 
\hspace{-3em}\frac{dD}{dt} &= k_{dn} \frac{f(N)^6}{x_{dn}^6 + f(N)^6}-\mu_{d} D \label{eq:mod3}\\ 
\hspace{-3em}\frac{dC_{a}}{dt} &= s_{c}+k_{cn} \frac{f(N+k_{cnd}D)}{1+f(N+k_{cnd}D)}-\mu_{c}C_{a}\label{eq:mod4}+u_a(t)
\end{align}} }

\noindent Set 
\begin{eqnarray}
	R=f(k_{np}P+k_{nn}N+k_{nd}D), \hspace{2mm} 
	f(x)=\frac{x}{1+(\frac{C_a}{c_\infty})^2} \nonumber 
\end{eqnarray}
Table 1 gives the reference parameter values. Note that the state variables $P(t)$, $N(t)$, $D(t)$, $C_a(t)$ and the control variables $u_p (t)$, $u_a (t)$ take nonnegative values $\forall t$.

\begin{figure}[!ht] 
\begin{center}
    \includegraphics[scale=0.47]{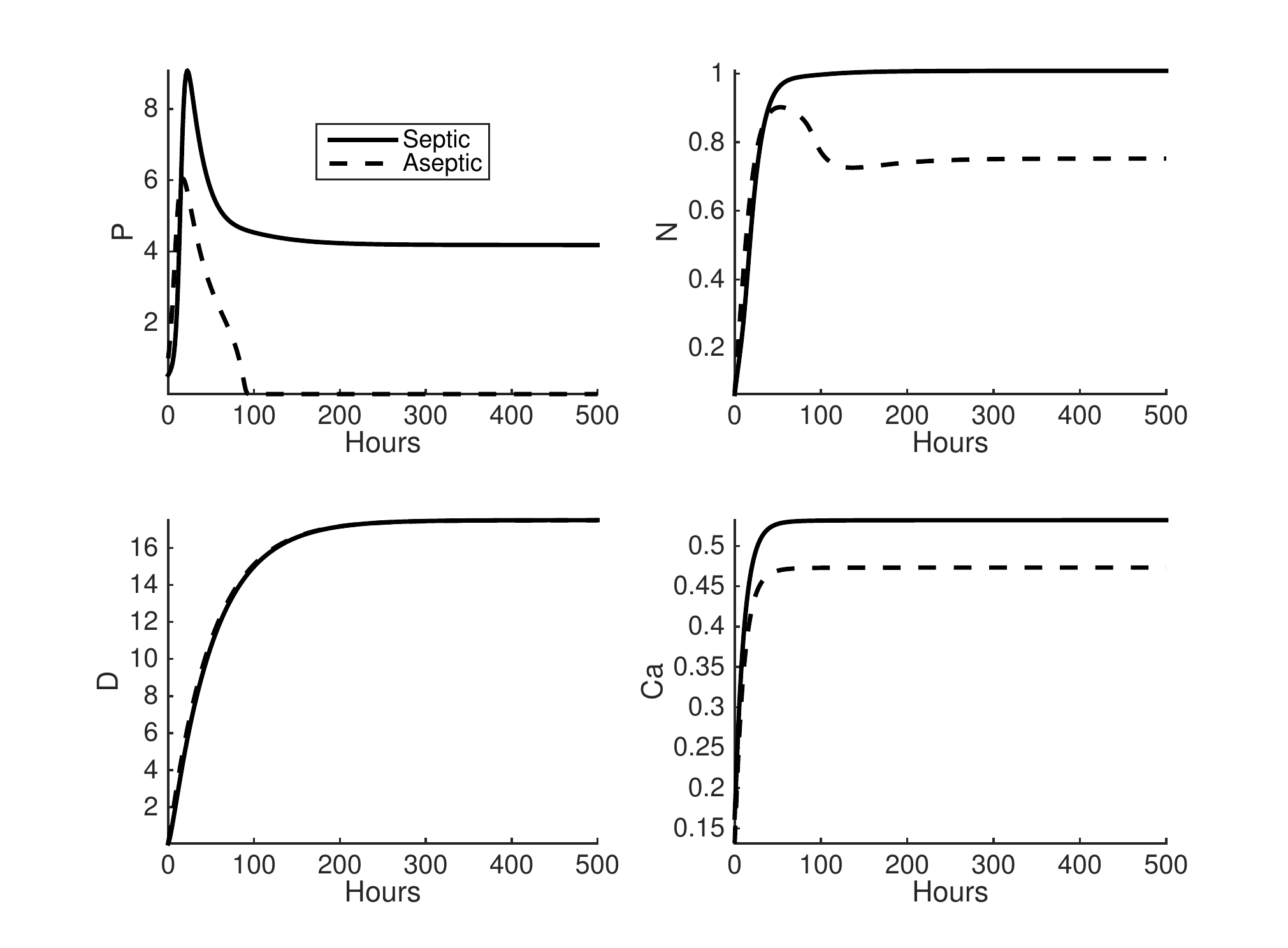} 
  \end{center}
  \caption{Natural (open loop) response  for patients 1 (septic) and patient 2 (aseptic) 
  }
  \label{Fig:fig1}
\end{figure}

\begin{itemize}
\item Equation (\ref{eq:mod1}) represents the evolution of the bacterial pathogen population $P$ that causes the inflammation. 
\item Equation (\ref{eq:mod2}) governs the dynamics of the concentration of a collection of early pro-inflammatory mediators such as activated phagocytes and the pro-inflammatory cytokines. They produce $N$. 
\item Equation (\ref{eq:mod3}) corresponds to tissue damage (D), which helps to verify the response outcomes. 
\item Equation (\ref{eq:mod4}) describes the evolution of the concentration of a collection of anti-inflammatory mediators $C_a$. 
\item See Tables I and II for the numerical values of the parameters and of the initial conditions.
\end{itemize}
The above model possesses three steady states: 
\begin{itemize}
\item one which corresponds to the healthy equilibrium, 
\item two which are associated respectively with a septic state and an aseptic one. 
\end{itemize}
Those properties agree with clinical observations: 
\begin{itemize}
\item The healthy equilibrium corresponds to $P=N=D=0$ and $C_a$ at a background level.  
\item A septic equilibrium is related to the situation where all mediators, $N$, $C_a$, and $D$ together with the pathogen $P$ are rather high. 
\item The patient is in an aseptic equilibrium when the values of $N$, $C_a$, $D$ are important, while the pathogen has been eliminated, \textit{i.e.}, $P=0$.  
\end{itemize}
See in Figure \ref{Fig:fig1} the results of two virtual patients with different initial conditions. The presence of  pathogen in the body stimulates inherently the activation of phagocytes (pro-inflammatory mediator). The resulting damage  is affected by the degree of inflammation which tries to eliminate the actual pathogen as quickly as possible. Note that  the actual anti-inflammatory mediator (cortisol and interleukin-10) can mitigate the inflammation and its harmful effect. The resting value $C_a$ is $0.125$ for the reference virtual patient. The patient is healthy when $D=0$ and $P=0$. He/she is considered to be dead when $D \geq 17$. When starting, \textit{e.g.}, from $[\begin{matrix}0.3& 0.0& 0.0& 0.0125 \end{matrix}]$, and allowing the pathogen to rise from a level of $P=0.3$ to  $P=0.6$, at some point the immune system is not strong enough to cope with the pathogen attack which will inevitably  attract the virtual patient to a septic or aseptic state (see Figure \ref{Fig:fig1}). Some intervention to stabilize the patient to its healthy equilibrium, \textit{i.e.}, to {\em homeostasis}, becomes mandatory. 

\begin{table}[!t]

  \begin{center}
    \includegraphics[scale=0.7]{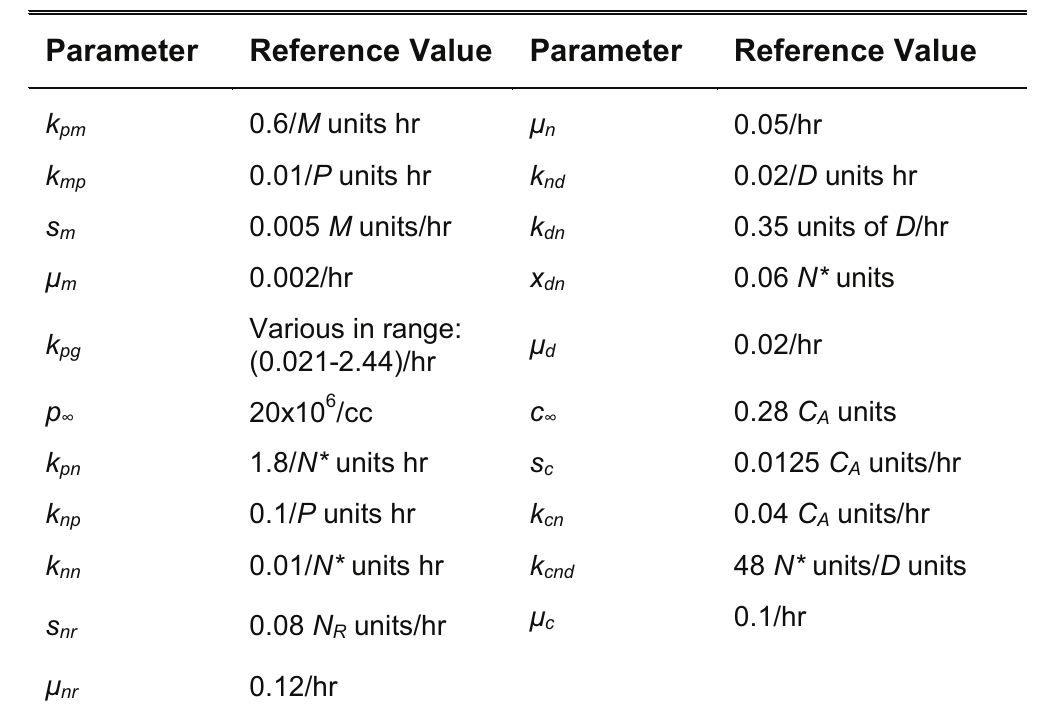}
    \caption{Reference parameters for the system (1)-(4)}
  \end{center}
  \label{table:tab1}
  \vspace{-3em}
\end{table}

\begin{table}[h!]
\centering
\begin{tabular}{cc}
\toprule
\bf{Parameter}    & \bf{ Parameter Ranges}  \\
\midrule 
$P_0$        & \,\,\,\,\,   0.0 -- 1.0      \\
$C_{A0}$     & \,\,\,\,\,   0.0938 -- 0.1563 \\
$k_{pg}$     & \,\,\,\,\,   0.3 -- 0.6       \\
$k_{cn}$     & \,\,\,\,\,   0.03 -- 0.05           \\
$k_{nd}$     & \,\,\,\,\,   0.015 -- 0.025       \\
$k_{np}$     & \,\,\,\,\,   0.075 -- 0.125      \\
$k_{cnd}$    & \,\,\,\,\,   36.0 -- 60.0       \\   
$k_{nn}$     & \,\,\,\,\,   0.0075 -- 0.0125  \\
\bottomrule 
\end{tabular} 
\caption{Variability of the model parameters}  \label{table:tab2}
\end{table}

\section{Model-free control\protect\footnote{See \cite{csm} for more details.}}\label{mfc}
\subsection{The ultra-local model}
Replace the unknown global description by the \emph{ultra-local model}:
\begin{equation}
\boxed{\dot{y} = F + \alpha u} \label{1}
\end{equation}
where
\begin{itemize}
\item the control and output variables are $u$ and $y$,
\item the derivation order of $y$ is $1$ like in most concrete situations,
\item $\alpha \in \mathbb{R}$ is chosen by the practitioner such that $\alpha u$ and
$\dot{y}$ are of the same magnitude.
\end{itemize}
The following explanations on $F$ might be useful: 
\begin{itemize}
\item $F$ is estimated via the measure of $u$ and $y$,
\item $F$ subsumes not only the unknown system structure but also
any perturbation.
\end{itemize}
\begin{remark}
In Equation \eqref{1} $\dot y$ is seldom replaced by $\ddot y$ (see, \textit{e.g.}, \cite{csm,ecc15}, and the references therein). Higher order derivatives were never utilized until today.
\end{remark}
\subsection{Intelligent controllers}
The loop is closed by an \emph{intelligent proportional controller}, or \emph{iP},
\begin{equation}\label{ip}
\boxed{u = - \frac{F - \dot{y}^\ast + K_P e}{\alpha}}
\end{equation}
where
\begin{itemize}
\item $y^\star$ is the reference trajectory,
\item $e = y - y^\star$ is the tracking error,
\item $K_P$ is the usual tuning gain.
\end{itemize}
Combining Equations \eqref{1} and \eqref{ip} yields:
$$
\dot{e} + K_P e = 0
$$
where $F$ does not appear anymore. The tuning of $K_P$, in order to insure local stability, becomes therefore quite straightforward. This is a major benefit when
compared to the tuning of ``classic'' PIDs (see, \textit{e.g.},
\cite{astrom,murray}, and the references therein), which
\begin{itemize}
\item necessitate a ``fine'' tuning in order to deal with the poorly known parts of the plant,
\item exhibit a poor robustness with respect to ``strong'' perturbations and/or system alterations.
\end{itemize}




\begin{figure*}
\begin{center}
\subfigure[Time evolution of $C_a$ and its reference trajectory (- -)]{
\resizebox*{8.3cm}{!}{\includegraphics{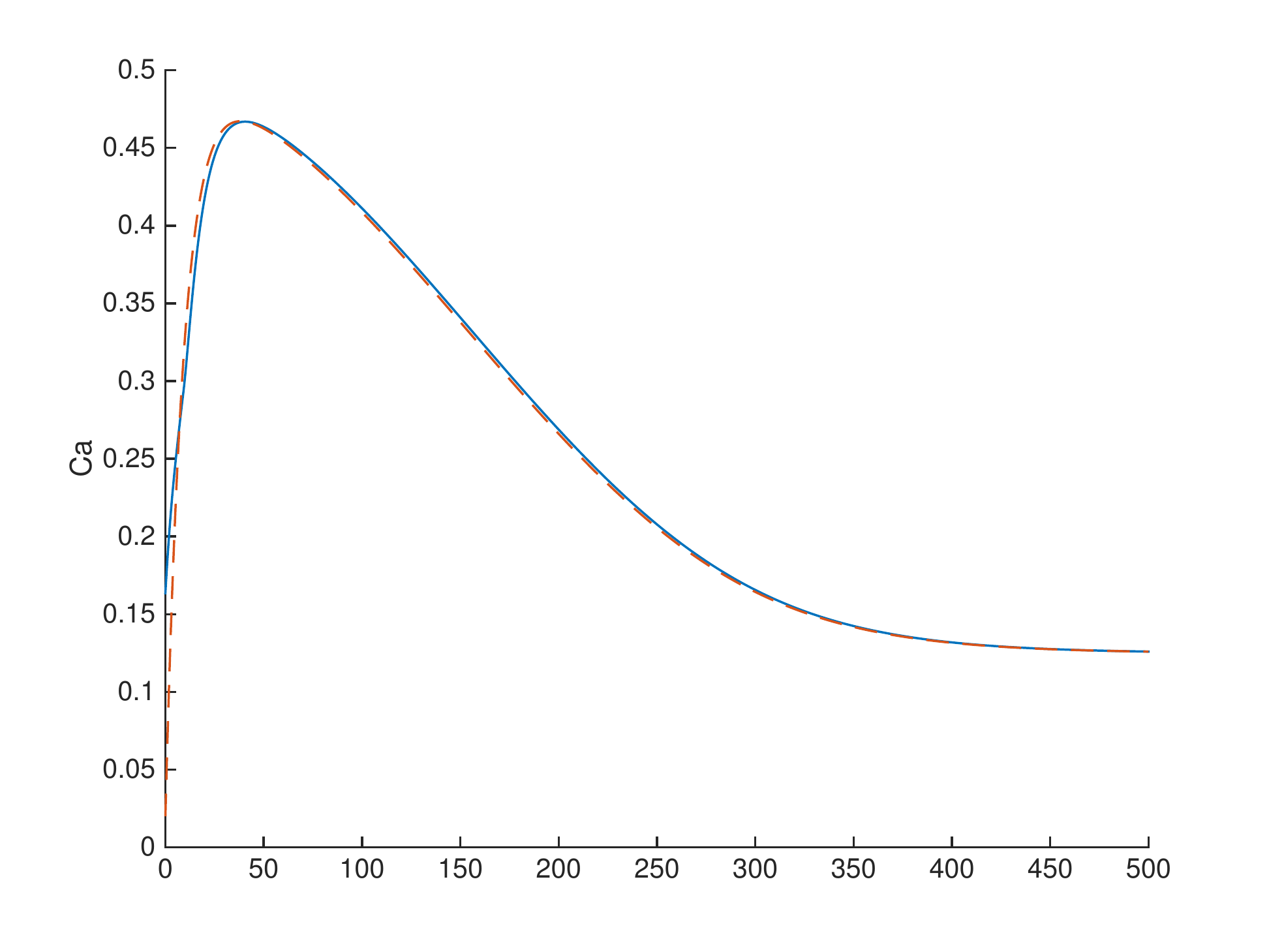}}}%
\subfigure[The evolution of $N$ and its nominal reference trajectory (- -)]{
\resizebox*{8.3cm}{!}{\includegraphics{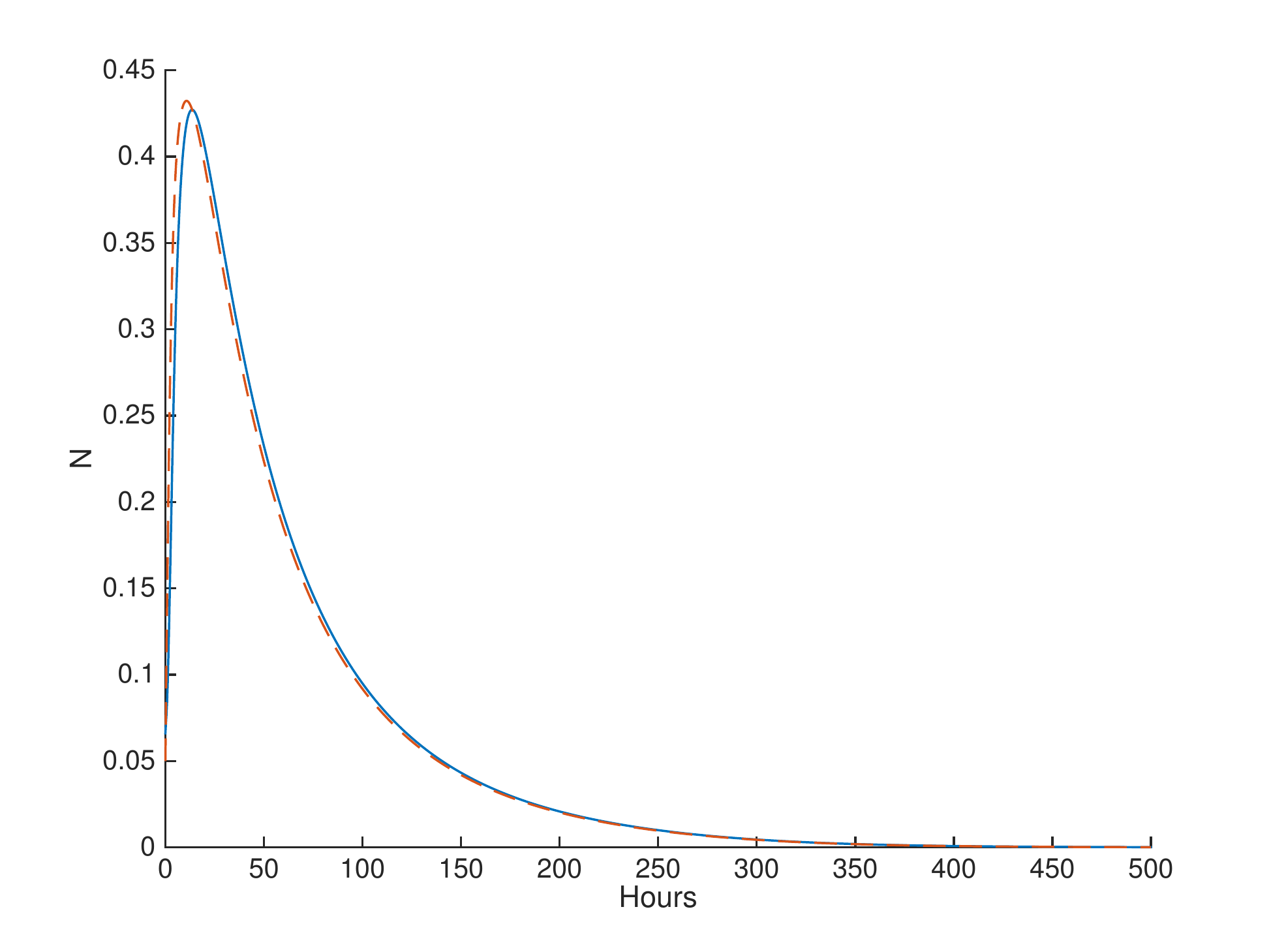}}}\\%
\caption{Reference trajectories $N$ and $C_a$ for both patients together with their closed loop response}%
\label{Fig:fig2}
\end{center}
\end{figure*}

\subsection{Estimation of $F$}\label{F}
The calculations below stem from new estimation techniques (see \cite{sira1,sira2}, and \cite{sira}).

\subsubsection{First approach}
The term $F$ in Equation \eqref{1} may be assumed to be ``well'' approximated by a piecewise constant function $F_{\text{est}} $. Rewrite then Equation \eqref{1}  in the operational domain (see, \textit{e.g.}, \cite{yosida}): 
$$
sY = \frac{\Phi}{s}+\alpha U +y(0)
$$
where $\Phi$ is a constant. We get rid of the initial condition $y(0)$ by multiplying both sides on the left by $\frac{d}{ds}$:
$$
Y + s\frac{dY}{ds}=-\frac{\Phi}{s^2}+\alpha \frac{dU}{ds}
$$
Noise attenuation is achieved by multiplying both sides on the left by $s^{-2}$. It yields in the time domain the realtime estimate, thanks to the equivalence between $\frac{d}{ds}$ and the multiplication by $-t$,
\begin{equation}\label{integral1}
{\small \boxed{F_{\text{est}}(t)  =-\frac{6}{\tau^3}\int_{t-\tau}^t \left\lbrack (\tau -2\sigma)y(\sigma)+\alpha\sigma(\tau -\sigma)u(\sigma) \right\rbrack d\sigma} }
\end{equation}

\subsubsection{Second approach}\label{2e}
Close the loop with the iP \eqref{ip}:
\begin{equation}\label{integral2}
\boxed{F_{\text{est}}(t) = \frac{1}{\tau}\left[\int_{t - \tau}^{t}\left(\dot{y}^{\star}-\alpha u
- K_P e \right) d\sigma \right] }
\end{equation}
\begin{remark}
Note the following facts: 
\begin{itemize}
\item integrals \eqref{integral1} and \eqref{integral2} are low pass filters,
\item $\tau > 0$ might be quite small,
\item the integrals may of course be replaced in practice by classic digital filters.
\end{itemize}
\end{remark}

\begin{remark}
A hardware implementation of the above computations is easy \cite{nice}.
\end{remark}

\section{Computer simulations}\label{virtual}
\subsection{Control design}
The state component $N$ (resp. $C_a$) in Equation \eqref{eq:mod2} (resp. \eqref{eq:mod4}) is
\begin{itemize}
\item easily measured, whereas it is difficult today to do it with $P$ and $D$ in Equations \eqref{eq:mod1} and \eqref{eq:mod4}.
\item mostly influenced by the control variable $u_p$ (resp. $u_a$).
\end{itemize}
Introduce therefore the  two Equations of type \eqref{1}:
{\begin{align}
\hspace{-3em}\dot{N} &= F_1 + \alpha_p u_p(t) \label{ump}\\ 
\hspace{-3em}\dot{C}_{a} &=  F_2 + \alpha_a u_a(t) \label{uma}
\end{align}} 

\noindent Let us emphasize that, like in \cite{toulon}, those two ultra-local systems may be ``decoupled'': they are considered as monovariable systems.\footnote{It should be nevertheless clear from a purely mathematical standpoint that $F_1$ (resp. $F_2$) 
is not necessarily independent of $u_a$ (resp. $u_p$).} The two corresponding iPs \eqref{ip} read
{\begin{align}
\hspace{-3em}u_p &= - \frac{F_1 - \dot{N}^\ast + K_{P1} e_p}{\alpha_p} \label{ipp}\\ 
\hspace{-3em}u_a &= - \frac{F_2 - \dot{C}_{a}^\ast + K_{P2} e_a}{\alpha_a} \label{ipa}
\end{align}} 
where the tracking errors are defined by
$$
e_p = N - N^\ast  \quad \text{and} \quad e_a = C_a - C_{a}^\ast
$$

\begin{figure*}[tbh!] 
\centering
\hspace{-2mm}
\includegraphics[width=5.01in]{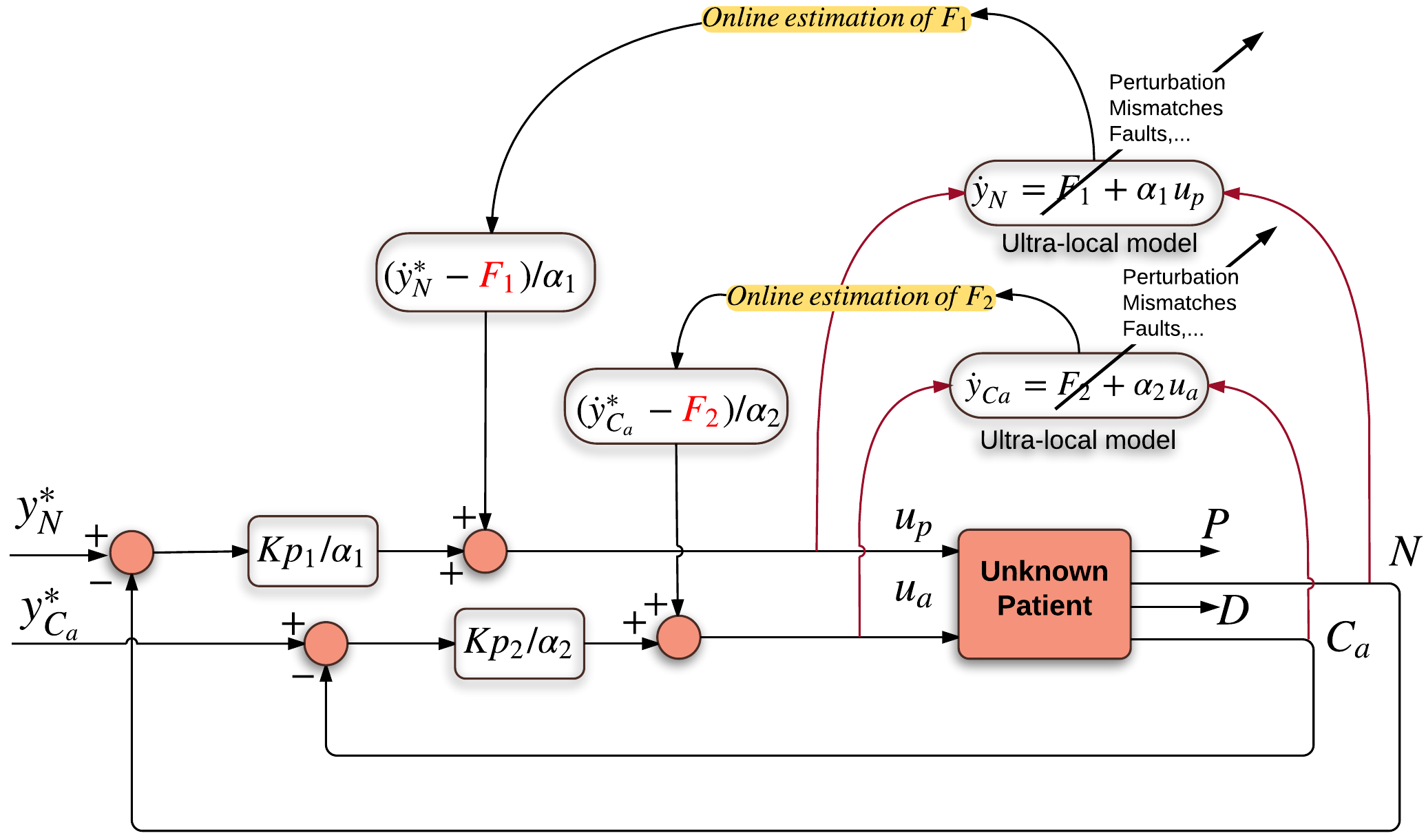}
\label{Fig:fig3}
\caption{Block diagram of our model-free control design.\label{scMFC}}
\end{figure*}

\noindent $F_1$, $F_2$ are estimated according to Section \ref{F}.
See  Figure \ref{scMFC} for the corresponding block diagram.
\subsection{Reference trajectories and results}
Two virtual patients are considered, the first (resp. second) one with a septic (resp. aseptic) outcome. The fact that a virtual patient may, or may not, return to an healthy state depends on the parameters and initial conditions. Their numerical
characteristics are given below:
\begin{enumerate}
\item Patient 1 (septic).-
\begin{itemize}
\item Initial conditions    $P(0) = 0.47360$, $N(0) = 0.0660$, $D(0) = 0.0477$, $C_a(0) = 0.1635$. 
\item Model coefficients $k_{pg} = 0.47846$, $k_{cn} = 0.0409$, $k_{nd} = 0.0242$, $k_{np} = 0.1211$, $k_{cnd} = 49.1243$, $k_{nn} =  0.012$.
\end{itemize}
\item Patient 2 (aseptic).-
\begin{itemize}
\item Initial conditions    $P(0) = 1.0017$, $N(0) = 0.0711$, $D(0) = 0.0732$, $C_a(0) = 0.1314$.
\item Model coefficients $k_{pg} = 0.4746$, $k_{cn} = 0.0386$, $k_{nd} = 0.0223$, $k_{np} = 0.1116$, $k_{cnd} = 46.3367$, $k_{nn} = 0.0112$.
\end{itemize}
\end{enumerate}
The reference trajectories of $N$ and $C_a$ are adjusted from Table I:
$$
N^\star=N_{\text{free}}.C_1,    \quad C_a^\star=({C_a}_{\text{free}}-0.125).C_2+0.125
$$
where 
\begin{itemize}
\item $N_{\text{free}}$ and ${C_a}_{\text{free}}$ correspond to the free trajectories of $N$ and $C_a$ for a healthy virtual patient,
\item $C_1$ and $C_2$ are suitable constants.
\end{itemize} 
 \begin{figure*}
\begin{center}
\subfigure[Time evolution of $P$ and $D$ for Patient 1 (--) and Patient 2 (- -)]{
\resizebox*{8.3cm}{!}{\includegraphics{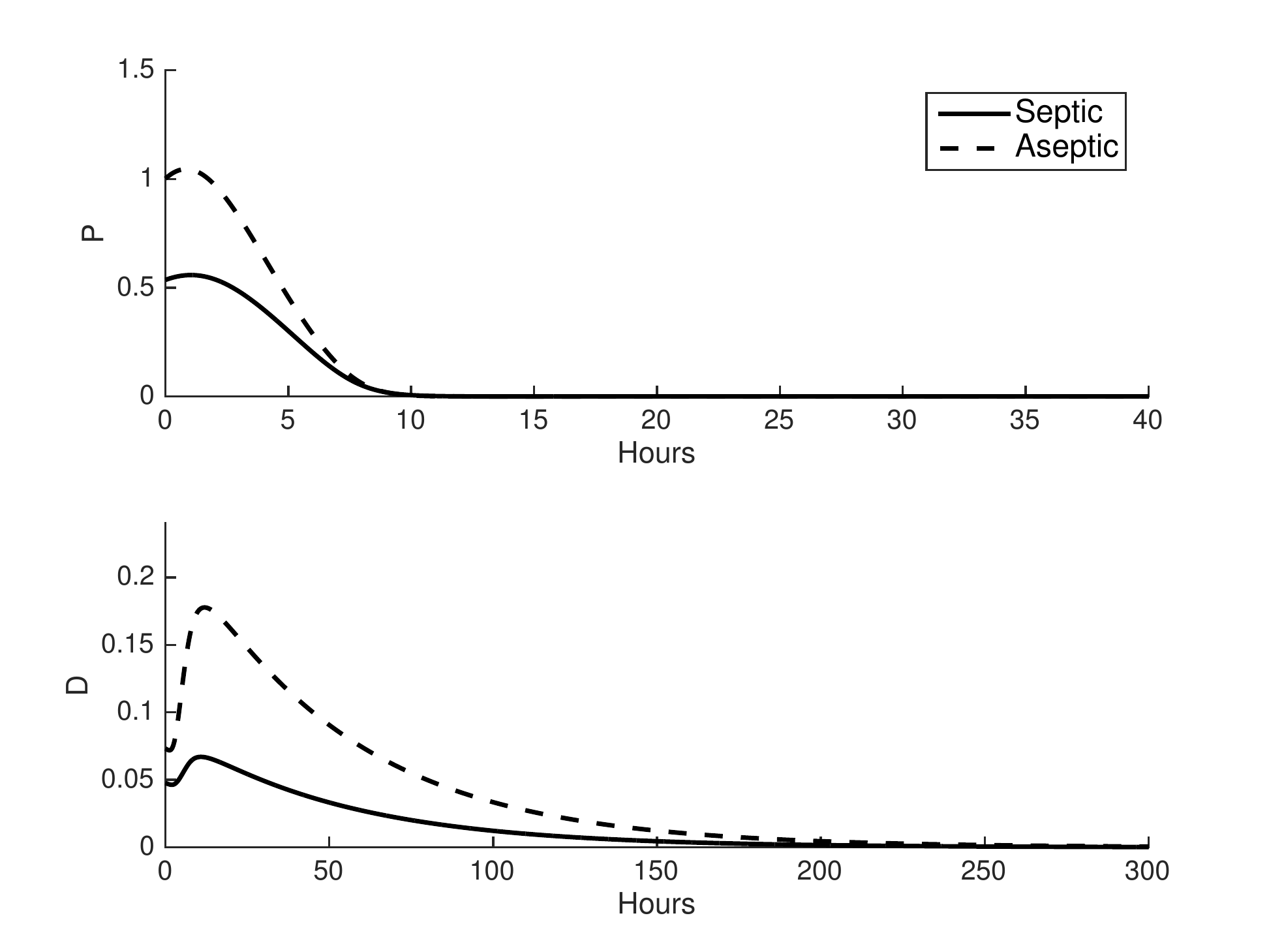}}}%
\subfigure[Control input $U_p$ and $U_a$  for septic and aseptic]{
\resizebox*{8.3cm}{!}{\includegraphics{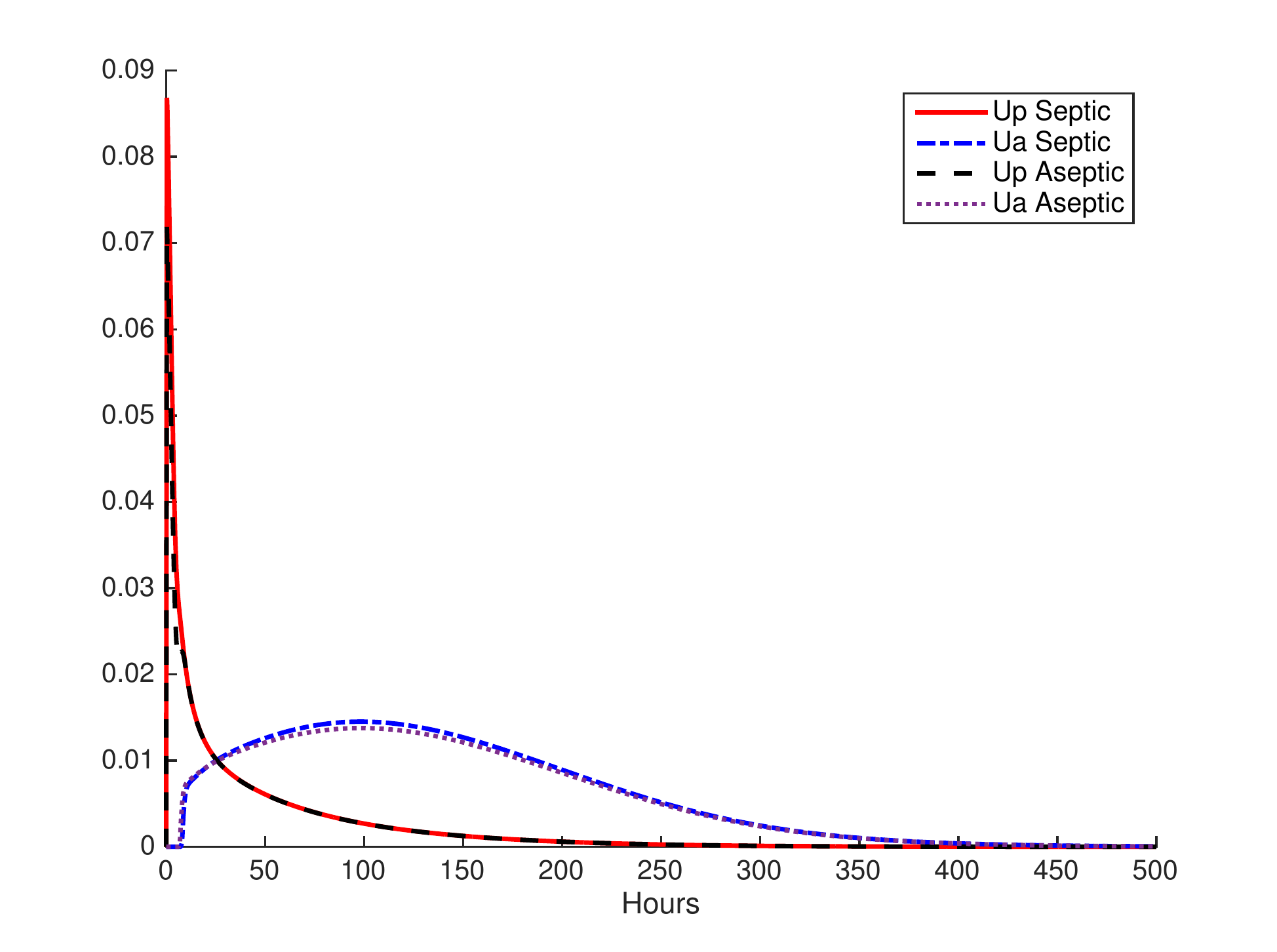}}}\\%
\caption{Patients 1 and 2 with therapy}%
\label{Fig:fig4}
\end{center}
\end{figure*}
We decided in our  scenario to amplify the trajectory corresponding to the concentration of the pro-inflammatory mediator. Therefore
$$
C_1 = 4, \quad C_2 = 1
$$
There are of course other possibilities for the reference trajectories. We could select higher amplitudes in order to heal most patients. The price to pay would be more drug injection and, therefore, more tissue damage. The simulation were performed 
\begin{itemize}
\item with a sampling time of $1$ minute, 
\item with $\alpha_p = \alpha_a = 2$ in Equations \eqref{ump}-\eqref{uma},
\item with $K_{P1} = K_{P2} = 0.47$ in Equations \eqref{ipp}-\eqref{ipa},
\item during $500$ hours.\footnote{Let us stress that our control objective was reached in less than $250$ hours.}
\end{itemize}
Figure \ref{Fig:fig4}(a)  shows clearly that we have been able to eliminate the pathogen and reduce the damage to zero using the generated doses displayed on the right hand side.
Many simulations show a quick rise in the 
pro-inflammatory mediator $N$. According to Figure \ref{Fig:fig2}(b), its maximum is reached after about 10 to 15 hours and is followed by an exponential decrease to zero. As shown by Figure \ref{Fig:fig2}(a), the analogous behavior of the anti-inflammatory mediator $C_a$ is much slower. Similar facts are observed with all patients who do not necessitate any treatment. The motivation for the choice of the reference trajectories should now become clear. 


The similarities of the generated doses can be partly explained by the same choice of the reference trajectory. In this case, it was enough to stabilize both patients. Observe that for each dose associated with an increase of the pro-inflammatory mediator a lower dose of anti-inflammation follows (see also \cite{bara2015,day2010using}). It may be explained by the fact  the immune system needs an initial boost of activated phagocytes in order to eliminate the pathogen threat. The resulting inflammation causes an increase of tissue damage, observed in Figure \ref{Fig:fig4}(a), which decreases after to zero thanks in part to the anti-inflammatory dose that is applied with a longer duration. Notice that injecting a larger dose of $U_a$ at the wrong time and with an inappropriate amplitude may foster the development of pathogen $P$ at rates that can drive the patient to a no-return point.

\balance

\section{Conclusion}\label{conclusion}

Our results should of course be further tested and developed. Future publications will emphasize 
\begin{itemize}
\item the robustness of our setting with respect to parameter variation and different initial conditions,
\item a deeper understanding of the choice of ``good'' reference trajectories,
\item the applicability of our approach to most types of inflammations and virtual patients.
\end{itemize}
The past success of model-free control in so many concrete situations should certainly be viewed as encouraging.

\end{document}